\title{Generalised Fermat Hypermaps and Galois Orbits}

\author{Antoine D. Coste\footnote{also 53 rue des meuniers, F--92150 Suresnes}\\
CNRS, UMR 8627, Building 210\\
Laboratory of Theoretical Physics\\
F--91405 Orsay Cedex, France\\
{\tt antoine.coste@m4x.org}
\and Gareth A. Jones\\
School of Mathematics\\
University of Southampton\\
Southampton SO17  1BJ,\\
U.K.\\
{\tt G.A.Jones@maths.soton.ac.uk}
\and  Manfred Streit\\
Oppenheimer Str. 36\\
D--60594 Frankfurt a.M.,\\
Germany\\
{\tt Manfred.Streit@bahn.de}
\and       J\"urgen Wolfart \\
Math. Sem. der Univ.\\ Postfach 111932\\
D--60054 Frankfurt a.M., Germany\\
{\tt wolfart@math.uni-frankfurt.de}}

\documentclass[12pt]{article}
\usepackage[latin1]{inputenc}
\usepackage{a4wide}
\usepackage{latexsym}
\usepackage{amsmath}
\usepackage{amssymb}
\usepackage{amsfonts}
\newtheorem{thm}{Theorem}
\newtheorem{lemma}{Lemma}
\date{}
\begin{document} 

\maketitle

\begin{abstract}
We consider families of quasiplatonic Riemann surfaces characterised
by the fact that --- as in the case of Fermat curves of exponent $n$ ---
their underlying regular (Walsh) hypermap is an embedding of the complete
bipartite graph $\,K_{n,n}\,$, where $\,n\,$ is an odd prime power. We show
that these surfaces, regarded as algebraic curves, are all defined over
abelian number fields. We determine their orbits under the action
of the absolute Galois group, their minimal fields of definition, and in
some easier cases their defining equations.
The paper relies on group-- and graph--theoretic results by G.~A.~Jones,
R.~Nedela and M.~\v{S}koviera about regular embeddings of the graphs
$K_{n,n}$ [JN\v{S}], and generalises the analogous results for maps
obtained in [JStW], partly using different methods.
\end{abstract} 



{\footnotesize Running head: Generalised Fermat hypermaps\\
Address for correspondence: J. Wolfart, Phone +49 (0)69 798 23423, Fax
+49 (0)69 798 28444,\\ wolfart@math.uni-frankfurt.de}

\section{Definitions and main results}

Riemann surfaces $X$ uniformised by subgroups $\Gamma$ of triangle groups 
$\Delta$ play a special role as {\em Bely\u{\i} surfaces}, that is, surfaces having 
a {\em Bely\u{\i} function}
$$\,\beta:X\to{\bf P}^1({\bf C})$$
ramified over at most three points in the Riemann sphere $\,{\bf P}^1({\bf C})\,$,
corresponding to the covering map
$$ \beta \; : \; \Gamma \backslash {\bf H} 
\; \to \; \Delta \backslash {\bf H} \;$$
where ${\bf H}$ is the hyperbolic plane. 
As first observed by Bely\u{\i} [Bel], the existence of such a function is 
equivalent to the property that $X$ --- as a smooth projective algebraic 
curve --- can be defined over a number field. Starting with Grothendieck's theory of 
{\em dessins d'enfants} [G], many interesting reformulations of Bely\u{\i}'s theorem
have been found,  see for instance [VS], [CIW], [JS], the recent survey in [W2], or the 
introduction in [JStW]. For the present paper the most important 
aspects are on the one hand the uniformisation by subgroups of triangle 
groups, and on the other hand the motivation from the theory of hypermaps
and their {\em Walsh representations}. The Walsh map $W(H)$ of a hypermap $H$ is a
bipartite graph embedded in a compact orientable surface, dividing it into
simply connected cells; the black and white vertices represent the hypervertices and
hyperedges of $H$, the edges represent incidences between them, and the cells
represent the hyperfaces. Every Bely\u{\i} function $\beta$ induces such a bipartite map
on a Riemann surface $X$: if we normalise its critical values to be $\,0, 1$ and
$\infty\,$ then $\,\beta^{-1}(0)\,$ and $\,\beta^{-1}(1)\,$ are the sets of white and
black vertices, and the connected components of the preimage of the real interval
$\,]0,1[\,$ are the edges of the graph. Conversely, every bipartite map (equivalently,
every hypermap) on a compact orientable surface arises in this way from a unique
holomorphic structure and a unique Bely\u{\i} function on the surface. 

In this situation, an important problem is that of relating the combinatorial properties
of the hypermap $H$ to the algebraic properties of the curve $X$, such as its
moduli field, Galois orbit, defining equations, etc. In general, this problem
is very difficult, but it is a little easier if the Bely\u{\i} function
is a regular covering, that is, $\beta$ is the quotient map $X\to G\backslash X$ by a
group
$G$ of holomorphic automorphisms of the Riemann surface $X$; this is equivalent to
$\Gamma$ being a normal subgroup of the triangle group $\Delta$, with
$G\cong\Delta/\Gamma$, and also to the hypermap $H$ being regular, that is, having an
automorphism group, isomorphic to $G$, acting transitively on the edges of the Walsh
map $W(H)$. Such surfaces $X$, known as {\em quasiplatonic surfaces}, have many
interesting properties,  see for instance [W2, Thm.~4]. In particular $G$ can be
identified with the Galois group of the extension of function fields corresponding to
$\beta$, and the Galois correspondence allows information about $G$ 
and its action on $X$ to be
translated into information about this extension (see [JStW], [St], [StW1] and [StW2]
for examples of this).

Here we consider these problems for a family of Riemann surfaces
$X=X(f;u,v,w)$, defined later in this section, which can be regarded as generalisations
of the well-known Fermat curves. Our main results are presented in Theorems~1, 2 and 3,
but before defining these surfaces and stating our results we will give some background
information to motivate our choice of these examples.

A simple and classic example of regularity is the Fermat curve of exponent $n$,
given in projective coordinates by $x^n+y^n+z^n=0$, with Bely\u{\i} function
$\beta([x,y,z])=(x/z)^n$. Here the bipartite graph is as symmetric as it could be,
namely the complete bipartite graph $K_{n,n}$; the group $G$ is the direct product of
two cyclic groups of order $n$, consisting of the automorphisms of $X$ which multiply
$x$ and $y$ by a pair of $n$th roots of unity, while $\Delta$ is the triangle group
$[n,n,n]$ and $\Gamma$ is its commutator subgroup.

In recent years, considerable progress has been made towards understanding all the
regular embeddings of complete bipartite graphs. Here we must distinguish
between the regularity of the hypermap $H$, which requires a group of automorphisms to
act transitively on the edges of $W(H)$, and the stronger condition of regularity of the
map $W(H)$, which requires it to act transitively on the {\em directed} edges of
$W(H)$, so that there is an additional automorphism reversing an edge and hence
transposing the vertex-colours. In the case where $n$ is an odd prime power, the regular
maps which embed the graph $K_{n,n}$ have been classified in [JN\v S], and the Galois
theory associated with these maps has been investigated in [JStW]. However, in the
context of dessins d'enfants, it is regularity of the hypermap, rather than the map,
which is the more interesting property. In fact, the methods used in [JN\v S] implicitly
determine the wider class of edge-transitive complete bipartite maps, or equivalently,
the regular hypermaps with Walsh map $K_{n,n}$, again for odd prime powers $n$. The aim
of this paper is to extend the results obtained in [JStW] by studying the action of the
absolute Galois group on the larger family of curves associated with these hypermaps.

It is shown in [JN\v S] that if a regular map is an embedding of $K_{n,n}$, for any $n$,
then the group $G$ of automorphisms preserving the orientation and the vertex-colours
has elements $x$ and $y$ (rotations around a black and a white vertex) such that

\begin{enumerate}

\item[(i)] $x$ and $y$ have order $n$,

\item[(ii)] $G=\langle x\rangle\langle y\rangle$ and $\langle x\rangle\cap\langle
y\rangle=1$, and

\item[(iii)] $x$ and $y$ are transposed by an automorphism of $G$.

\end{enumerate}

Conversely, every group $G$ with such a pair $x, y$ arises in this way: the
black and white vertices can be identified with the cosets of $\langle x\rangle$ and
$\langle y\rangle$, and the incident edges with the elements of those cosets,
cyclically ordered by successive powers of $x$ or $y$. Isomorphism classes of maps
correspond to orbits of ${\rm Aut}\,G$ on pairs $x, y$ satisfying these conditions. If
condition (iii) is omitted, we obtain a similar group-theoretic characterisation of the
edge-transitive embeddings of $K_{n,n}$, or equivalently the regular hypermaps
associated with this graph.

In the case where $n$ is an odd prime power $p^e$, it has been shown in [JN\v S] that
the regular embeddings of $K_{n,n}$ correspond to the groups
$$G = G_f \; := \; \langle g,h\mid g^n=h^n=1,\, g^{-1}hg=h^{1+p^f}\rangle $$
where $\,f=1,2,\ldots, e\,$. Note that each such group $G_f$ is a semidirect product $\,C_n
\rtimes C_n\,$ of a cyclic normal subgroup $\langle h\rangle$ by a cyclic subgroup
$\langle g\rangle$, that different values of $f$ give non-isomorphic groups, and that
$G_e$ is the direct product  $\,C_n \times C_n\,$. Now the arguments used in Sections
3--6 of [JN\v S] to show that $G=G_f$ for some $f$ depend only on conditions (i) and
(ii), and not condition (iii), so in fact they show that for odd prime powers $n$ the
edge-transitive embeddings of $K_{n,n}$ are also associated with these groups $G_f$: the
only difference is that in this case $x$ and $y$ need not be transposed by an
automorphism. As shown in [JN\v S], elements $x$ and $y$ of $G_f$ satisfy conditions
(i) and (ii) if and only if they generate $G_f$, or equivalently
$$x=g^uh^s \; , \quad  y=g^vh^t \quad {\rm with} \quad ut-vs \not\equiv 0 \bmod p\,\;,$$
in which case $x, y$ and $xy$ all have order $n$.

We need each of our Riemann surfaces $X$ to be uniformised by the
torsion-free kernel $\Gamma$ of an epimorphism $\,\theta:\Delta\to G_f\,$ for some
triangle group $\Delta$, so to allow this we take $\Delta$ to be the triangle group
$$ \Delta \,=\, [ n,n,n] \,:=\, 
\langle \gamma_0,\gamma_1,\gamma_{\infty} \,|\, \gamma_0^n=\gamma_1^n 
=\gamma_{\infty}^n=1=\gamma_0\gamma_1\gamma_{\infty} \rangle\,, $$
acting in the usual way on the upper half plane (or the complex plane if $n=3$). The
surface groups $\Gamma$ in question are defined to be the kernels  
$\Gamma_{u,v,w}$ of epimorphisms $\,\theta=\theta_{u,v,w}:\Delta \to G_f\,$ given by 
$$ \gamma_0 \mapsto x:=g^uh^{s_0} \; , \quad \gamma_1 \mapsto y:=g^vh^{s_1} 
\;, \quad \gamma_{\infty} \mapsto (xy)^{-1}=g^wh^{s_{\infty}} \;, $$ 
where the exponents of the images satisfy certain obvious conditions: 
since $G_f$ has exponent $n$ [JN\v S], the relation $\,\gamma_0 \gamma_1 \gamma_{\infty}
= 1\,$ shows that $\,\theta\,$ is a well-defined homomorphism if and only if 
$$ u+v+w  \; \equiv \; 0 \; \bmod n $$
and
$$\,s_0q^{v+w} +s_1q^w+s_{\infty} \equiv 0\; \bmod n$$
where $ \, q:= 1+p^f\,$. 
This homomorphism is surjective if and only if
$$\,us_1-vs_0 \not\equiv 0 \bmod p\,$$
[JN\v{S}, Prop.~12], in which case the images of the generators $\gamma_0, \gamma_1$
and $\gamma_{\infty}$ have order $n\,$, so the kernel is torsion-free and is therefore a
surface group. Given $u, v$ and $w$ satisfying $\, u+v+w  \; \equiv \; 0 \; \bmod n\,
$, one can choose $s_0, s_1$ and $s_{\infty}$ to satisfy the other two conditions
provided at least one of $u, v$ and $w$ is not divisible by $p$. We will see later that
the exponents  $\,s_0,s_1\,$ and $\,s_{\infty}\,$ play only a minor role in our
calculations, so we omit them in the notation for the kernels and the resulting
surfaces. 

Under these conditions, specifically $u+v+w \equiv 0 \bmod n\,,$ at
least one term not divisible by  $p$, we consider the Riemann surfaces 
$\,X(f;u,v,w) := \Gamma_{u,v,w} \backslash {\bf H}\,$ and prove

\begin{thm}
For $\,i=1,2\,$ let $\,u_i, v_i , w_i \in {\bf Z} \,$, not all 
divisible by $p\,$, satisfy the congruences $\,u_i+v_i+w_i \equiv 0
\bmod n\,$. Then 
the Riemann surfaces $X(f;u_1,v_1,w_1)$ and $X(f;u_2,v_2,w_2)$ are 
isomorphic if and only if $\,u_1,v_1,w_1\,$ are congruent 
$\, \bmod \; p^{e-f}\,$ to a permutation of $\,u_2,v_2,w_2\,$.
\end{thm}

These surfaces $X(f;u,v,w)$ form a complete list of all surfaces with 
{\em regular} hypermaps (or dessins) based on $K_{n,n}\,$, where {\em regular\/} now 
means that there is a group of holomorphic automorphisms of the surface 
(in our case $G_f\,$) acting transitively on the edges and preserving the vertex
colours.  In [JStW] the special case of regular {\em maps} 
is treated, i.e.~those with additional holomorphic involutions of the surfaces 
reversing the vertex colours; we will see that these cases form 
the special  subfamily defined by $\,u=v\,$. We call the dessins on the 
surfaces $X(f;u,v,w)$ considered here {\em generalised Fermat hypermaps} since in the
extremal case $\,f=e\,$ they live on $X(e;1,1,1)\,$, and this is the 
Fermat curve with exponent $\,n= p^e\,$. In fact, all curves
$X(f;u,v,w)$ have Weil curves (see (1) in Thm.~5) as common 
quotients with Fermat curves. Both Weil and Fermat curves and their
Jacobians have attracted recent interest among physicists, see [BCIR].

The aim of this paper is the same as in [JStW], namely to study the action of the
absolute Galois group, but now on this greatly extended family of curves. Two of the
most interesting points of dessin theory have been  already mentioned: 
\begin{enumerate}
\item[ ] as algebraic curves, Bely\u{\i} surfaces can be defined 
over number fields; 
\item[ ] they are uniquely defined by their underlying hypermap. 
\end{enumerate}
 
In theory, therefore, all interesting information about such a curve
should be encoded in the combinatorial or group-theoretical data 
given by the hypermap. In particular, it should be possible to use 
the hypermap to obtain defining equations (a process which might lead to extremely
technical  questions), or at least the (minimal) field of definition and the behaviour 
under algebraic conjugation, acting on the coefficients of the 
equations of the curve and of the Bely\u{\i} function: in other words, to determine
the Galois orbit of the dessin. To describe our main result, define 
$\,\eta := e^{2 \pi i/p^{e-f}} \,$ and recall that the Galois group 
$\,{\rm Gal \,}{\bf Q}(\eta)/{\bf Q}\,$ is isomorphic to the multiplicative 
group of units $\,({\bf Z}/p^{e-f}{\bf Z})^*\,$,
where the action of the residue class $[k]$ is given by
$\,\eta \mapsto \eta^{k}\,$. 

\begin{thm}
Under the hypotheses described above, two curves $X(f;u_1,v_1,w_1)$ 
and\\ 
$X(f;u_2,v_2,w_2)$ are Galois conjugate if and only if there is a 
$\,k \in {\bf Z}\,$ not divisible by $p$ such that $\,ku_1,kv_1,kw_1\,$ 
are congruent $\,\bmod \; p^{e-f} \,$ to a permutation of $\,u_2,v_2,w_2
\,$. The curves $X(f;u,v,w)$ can all be defined over subfields of ${\bf Q}(\eta) 
\;$.
\end{thm}

The precise determination of the minimal field of definition of the curves 
turns out to be rather more technical. Recall that the Galois group
always has a subgroup $\,\{ \pm 1 \}\,$ of order $2$ whose fixed field
under the Galois correspondence is the maximal real subfield 
$\,R := {\bf Q}(\cos 2\pi/p^{e-f}) \subset  {\bf Q}(\eta)\,$, 
and some more number theory shows that there is a subgroup $\,\{ 1,k,k^2 \} \,$ 
of order $3\,$, i.e.~with $\,k^3=1\,$ or better (to avoid the trivial 
solution $\,k=1\,$) $\,1+k+k^2 = 0\,$ if and only if $\,p \equiv 1 \bmod 3\,$. 
If we let  $K$ denote the subfield of ${\bf Q}(\eta)$ fixed by this group, we can 
give the field of definition in the following form.

\begin{thm}
Under the hypotheses given above, the (minimal) field of definition 
of the curve $X(f;u,v,w)$ is
\begin{enumerate}
\item $R \subset {\bf Q}(\eta)\,$ if and only if one of the parameters 
$\,u,v,w \equiv 0 \bmod p^{e-f}\,$,
\item $K \subset {\bf Q}(\eta)\,$ if and only if the parameters are of the 
form $\,u,ku,k^2u\,$ with some $u$ and $\,k \in {\bf Z}\,$, both coprime to  
$p\,$, such that $\,1+k+k^2 \equiv 0 \bmod p^{e-f}\,$, 
\item ${\bf Q}(\eta)\,$ in all other cases. 
\end{enumerate}
\end{thm}

Two further theorems treat the full automorphism groups of the curves 
(end of Sec.~2) and the explicit determination of their equations (Sec.~4). 
The main line of reasoning is similar to that of 
[JStW], but some new ideas were needed as well, in particular for 
the proofs of Theorems 1, 3 and 4, and in other cases we will give 
alternative proofs. 

The central point of the paper --- contained in 
Theorems 2 and 3 --- is the information that quasiplatonic curves 
with regular hypermaps based on the complete bipartite graphs $K_{n,n}$ 
are defined over abelian number fields. This result fits well into a more 
general framework since some other series of hypermaps with similar properties 
are known: in [StW1, Thm.~1] an analogous statement is proved for 
quasiplatonic curves whose hypermaps are based on complete bipartite 
graphs $K_{p,q}\,, \,q\,$ and $\,p \equiv 1 \bmod q\,$ primes $\,>3\,$ 
and $\,\not= 7\,$, and whose automorphism groups are also semidirect products 
$\,C_p \rtimes C_q\,$ of cyclic groups. Another series of examples are  
the hypermaps treated in [StW2] coming from bipartite graphs for finite 
cyclic projective planes of order $q$ with a Singer group $\,C_l,\, 
l=q^2-q+1\,$. These hypermaps are in general not regular, but in many 
cases (maybe with the only exception $\,q=5\,,\, l= 21\,$) they are 
uniform, i.e.~all vertices with valency $q\,$, all faces with valency $2l\,$, 
and with automorphism group containing $C_l\,$, see the discussion of 
the {\em Wada property} in [StW2, Sec.~3]. The graphs are not 
themselves bipartite  
but their dual graphs are complete bipartite graphs of type 
$K_{q,l}\,$, and the curves are again defined over abelian number fields 
[StW2, Prop.~6]. It is therefore a natural question whether sufficiently regular 
dessins based on complete bipartite graphs always lead to curves defined over 
abelian number fields.

{\small The first and the second author thank the Mathematisches 
Seminar of Frankfurt University 
for its scientific hospitality during the academic year 2005/06.}


\section{Isomorphisms and automorphisms}

To prove Theorem 1, recall that $\Delta$ is a normal subgroup of a
triangle group $\,\tilde{\Delta} = [2,3,2n]\,$ whose order $3$ generator
$\delta$ can be chosen so that its fixed point is the hyperbolic
midpoint of the fixed points $\,z_0,z_1,z_{\infty}\,$ of (respectively) $\,\gamma_0, \gamma_1,
\gamma_{\infty}\,$, see [Si]. Conjugation by $\delta$ provides a cyclic
permutation of the generators $\gamma_i$ of $\Delta$ conjugating
the kernel $\Gamma_{u,v,w}$ into $\Gamma_{v,w,u}\,$. In other words,
cyclic permutations of the parameters $\,u,v,w\,$ only induce isomorphisms
of the surfaces $X(f;u,v,w)\,$. 

We may therefore assume that at least the parameters $u$ and $v$ are
coprime to $p\,$. If we compose the epimorphism $\,\Delta \to G_f\,$
with an automorphism of $G_f\,$, the kernel does not change, therefore 
[JN\v{S}, Prop.~16] allows us to simplify the images of our generators
considerably. Since all automorphisms are given by $\,g \mapsto g^ih^j\,,\, h
\mapsto g^kh^l\,$ with $\,i \equiv 1 \bmod p^{e-f} \,,\; k \equiv 0 \bmod
p^{e-f}\,,\; l\,$ not divisible by $p\,$, we may always assume that $\,s_0=0\,,\, s_1 =
1\,$. Moreover, it is evident that the isomorphism class of $X(f;u,v,w)$
depends only on the residue classes of $\,u, \, v \bmod p^{e-f}\,$, and
since $w$ is uniquely determined by $u$ and $v\,$, we may consider 
$w$ also as a residue class $\; \bmod \; p^{e-f}\,$. 

Our epimorphism $\,\Delta \to G_f\,$ is now normalised to satisfy 
\[ \gamma_0 \mapsto g^u \; , \quad \gamma_1 \mapsto g^vh 
\;, \quad \gamma_{\infty} \mapsto g^wh^{s} \;,\quad u,v \not\equiv 0
\bmod p \;, \]
and the remaining exponent $s$ is determined by $u$ and $v\,$. This is
because 
$\,\gamma_0 \gamma_1 \gamma_{\infty} = 1\,$ and the defining relations
of $G_f$ (in particular $\,hg=gh^q\,$ with $\,q= 1+p^f\,$) imply
\[ g^{u+v}hg^wh^s = 1 \;, \quad \mbox{hence} \quad q^w+s \equiv 0 \bmod
n \;,\] 
where in fact $\,q^w \equiv q^{w'} \bmod n\,$ 
if and only if $\,w \equiv w' \bmod p^{e-f}\,$. 

To prove that all permutations of $\,u,v,w\,$ lead to isomorphic curves,
it is now sufficient to prove that $\,X(f;u,v,w) \cong X(f;v,u,w)\,$. As
in the initial argument concerning cyclic permutations there is also an
order $2$ generator $\delta$ of
$\tilde{\Delta}$ whose fixed point is the hyperbolic midpoint
$z_{\frac{1}{2}}$ of the
fixed points $\,z_0,z_1\,$ of $\,\gamma_0, \gamma_1 \in \Delta\,$. Conjugation by
$\delta$ transposes these two generators and sends $\gamma_{\infty}$ 
to $\,\gamma_0^{-1} \gamma_{\infty} \gamma_0\,$. In other words, $\delta$
conjugates the kernel $\Gamma_{u,v,w}$ to the kernel of the
epimorphism determined by
\[ \gamma_0 \mapsto g^vh \; , \quad \gamma_1 \mapsto g^u 
\;, \quad \gamma_{\infty} \mapsto g^{-u}g^wh^sg^u = g^w h^{sq^u}\;.\]
Using the relation $\,hg=gh^q\,$ and composition with a suitable
automorphism of $G_f$ a straightforward calculation shows that this
kernel coincides with $\Gamma_{v,u,w}\,$. 

Are there more isomorphisms between $X(f;u,v,w)$ and other curves of
this family? 
If so, their surface groups would be conjugate in ${\rm PSL}_2 {\bf
  R}\,$, and by [GW, Thm.~9] they would be conjugate even in the
normaliser $N(\Delta)$ 
(all normalisers taken in ${\rm PSL}_2 {\bf R}\,$) or in the normaliser
$\,N(\hat{\Delta})\,$ where $\hat{\Delta}$ denotes the normaliser 
of the surface
group $\Gamma_{u,v,w}\,$. All possibilities are well known by
Singerman's work [Si]. For the first possibility we know that $\,N(\Delta) =
\tilde{\Delta} = [2,3,2n]\,$, and we already know that conjugation by
elements of $\tilde{\Delta}$ only causes permutations of $\,u,v,w\,$  ---
up to congruences $\; \bmod \; p^{e-f}\,$. In the second case we 
obviously have $\,\Delta \subseteq \hat{\Delta}\,$, therefore $\,\hat{\Delta} 
\subseteq \tilde{\Delta}\,$ by [Si] again with one possible exception: the
triangle group $\,\Delta = [9,9,9]\,$ is contained with index $12$ in
the triangle group $\,[2,3,9]\,$. As a maximal triangle group, this 
group is its 
own normaliser, therefore all conjugations in question would again lead
only to
automorphisms of curves. This finishes the proof of Theorem 1. In 
Section 3, we will sketch a different argument for the {\em only if}
part. 

In the last mentioned example we can prove moreover that 
the group $\,[2,3,9]\,$ can never be the 
normaliser of one of our surface groups: if $\,(e=) f=2\,$ we have the 
Fermat curve of exponent $9$ where $\, 
N(\Gamma_{1,1,1}) = \hat{\Delta}= \tilde{\Delta}=[2,3,18]\,$, and if
$\,f=1\,$, 
\[ N(\Gamma_{u,v,w}) \; = \; \hat{\Delta} \;=\; [2,3,9] \]
would imply the existence of some cyclic automorphism of $G_f$ of order
$3$ with the effect 
\[ g^u \mapsto g^vh \mapsto g^wh^s \mapsto g^u \;, \]
see [BCC, case T6]. The existence of such an automorphism, together with  
our conditions on $\,u,v,w\,$ and the shape of automorphisms
of $G_f\,$, see again [JN\v{S}, Prop.~16], implies that $\,u \equiv v \equiv 
w \equiv 1 \,$ or $2 \bmod 3\,$. But then all permutations of $\,u,v,w\,$ 
induce automorphisms of the curve, induced by the action of 
$\,\tilde{\Delta}/\Delta\,$, therefore we have in fact 
\[ N(\Gamma_{u,v,w}) \;= \; \tilde{\Delta} \;=\; [2,3,18] \;.\] 

The same reasoning generalises to 
arbitrary prime powers $n$ where we can conclude similarly that the full
automorphism group of the curve always lifts to a triangle group
$\hat{\Delta}$ between $\Delta$ and $\tilde{\Delta}\,$. The automorphism
group is an
extension of $G_f$ by that subgroup of the permutation group $\,S_3 \cong
\tilde{\Delta}/\Delta\,$ leaving invariant the triple $\,(u,v,w) \bmod
p^{e-f}\,$. As an example, take the triples $\,(u,u,-2u)\,$ giving
surfaces with an additional holomorphic involution transposing the vertex
colours of the dessin, as discussed in [JStW], or the Fermat curves with
$\,e=f\,$ and triples $\,(1,1,1)\,$, or the special case $\,(u,u,u)\,$ if
$\,p=3\,,\; f=e-1\,$, see [JStW, Lemma 3b]. Summing up, we get the following 
generalisation of [JStW, Lemma 3]. 

\begin{thm}
The automorphism group of the curve $X(f;u,v,w)$ is an extension of 
$G_f$ by the permutation subgroup of $S_3$ leaving invariant the 
triple $\,(u,v,w) \bmod p^{e-f}\,$. 
\end{thm} 


\section{Galois action}

To prove Theorem 2, we start with the normalisation developed in Section 1,
i.e.~representing $X(f;u,v,w)$ for a fixed odd prime power 
$\,n= p^e > 3\,$ and fixed $\,f\, , \, 1\le f \le e\,$, as the quotient 
$\,\Gamma_{u,v,w} \backslash {\bf H}\,$ where $\Gamma_{u,v,w}$ denotes the 
kernel of the epimorphism $\,\Delta \to G_f \,$ determined by
\[ \gamma_0 \mapsto g^u \; , \quad \gamma_1 \mapsto g^vh 
\;, \quad \gamma_{\infty} \mapsto g^wh^{s} \;,\quad u+v+w \equiv 0 \bmod
n \; , \quad u,v \not\equiv 0 \bmod p \;, \]
Later we will have to compare the actions of the generating pairs $\,g,\,gh\,$
on all surfaces in question. To simplify this comparison, 
we change the epimorphism by composition with an automorphism 
of $G_f$ so that it is now given by
\[ \gamma_0 \mapsto g^u \;, \quad \gamma_1 \mapsto (gh)^v \; , \quad 
\gamma_{\infty} \mapsto g^wh^t \;.\]
This change can be justified in the same way as [JStW, Lemma 1]. 
Observe that the exponents 
$\,u,v,w\,$ remain unchanged and uniquely determine the exponent $t\,$. 

Any Galois conjugation sends $X(f;u,v,w)$ onto some 
other member $X(f;x,y,z)$ of the family because automorphism groups are sent to isomorphic 
automorphism groups and ramification orders stay invariant, so we have 
to look for Galois orbits only inside our families of curves with fixed 
$f\,$. The treatment of Galois action on this family relies on an 
idea first developed by Streit [St] for the case of Macbeath--Hurwitz 
curves, i.e.~the use of {\em multipliers}, used also for other families 
of Bely\u{\i} surfaces ([StW1], [StW2], [JStW]). 
To recall their definition, let $a$ be an automorphism of a Riemann 
surface $X$ with fixed point $P\,$. If $z$ is a local coordinate on 
$X$ in a neighbourhood of $P$ with $\,z(P)=0\,$ then 
\[ z \circ a \; = \; \xi z \,+\quad \mbox{higher order terms in} 
\quad z \;,\] 
and we call $\xi$ the {\em multiplier of $a$ at} $P\,$. Clearly, if $a$ is an
automorphism of order $n\,$ then $\xi$ 
is a $n$-th root of unity independent of the choice of $z$. Then (see 
[StW2, Lemma 4]) we have

\begin{lemma}
If $X$ has genus $\,g > 1\,$ and is defined over a number field, then
$a$ is also defined over a number field, $P$ is a 
$\overline{{\bf Q}}$-rational point of $X$, and for all 
$\,\sigma \in {\rm Gal}\,\overline{\bf Q}/{\bf Q}\,$, Galois 
conjugation of the coefficients by $\sigma$ gives an automorphism $a^{\sigma}$ 
of $X^{\sigma}$ with  multiplier 
$\sigma(\xi)$ at its fixed point $P^{\sigma}\,$. 
\end{lemma} 

A {\em proof} different from that given in 
[StW2] is the following. If $a$ were not 
defined over the field $\overline{{\bf Q}}\,$ there would be infinitely many different 
automorphisms $a^{\sigma}$ where $\sigma$ runs over field automorphisms 
of ${\bf C}$ fixing the field of definition of $X\,$, in contradiction to 
the finiteness of the automorphism group of $X\,$. Therefore the fixed 
point $P$ is also defined over $\overline{{\bf Q}}$, 
and there is a $\overline{{\bf Q}}$--rational function on $X$ 
unramified at $P$ serving as a local variable $z$ with $\,z(P)=0\,$. Then,
the function $z \circ a$ can be written as $\,\xi z + r\,$ where
$r$ is again a $\overline{{\bf Q}}$--rational function on $X$ with a
zero of order $\,>1\,$ at $P\,$, and this vanishing order is respected
by Galois conjugation. 

\begin{lemma}
In the normalisation chosen above, $g$ and $gh$ have $p^f$ fixed points 
on $X(f;u,v,w)\,$, respectively.  
\end{lemma}

Since $u$ and $v$ are coprime to $\,n=p^e\,$, it is sufficient to 
prove this claim for $g^u$ and $(gh)^v\,$ instead of $g$ and $gh\,$. For 
both, the arguments given in the proof of [JStW, Lemma 6]  
are valid without any change, so we omit the details. $\Box$ 

By the construction of the embedding of $K_{n,n}$ into $X(f;u,v,w)\,$, 
the fixed points of $g$ and of $gh$ consist of white and black vertices 
of the graph respectively, i.e.~under the canonical map $\,{\bf H} \to 
X(f;u,v,w)\,$ they come from certain points in the $\Delta$--orbit of 
the fixed points $\,z_0,z_1\,$ of $\gamma_0$ and $\gamma_1\,$. More precisely, they 
form orbits under the subgroup 
$\langle h^{p^{e-f}} \rangle$ which lies in the center of $G_f\,$. Therefore 
as in [JStW, Lemma 7] we get 
 
\begin{lemma}
Let $\zeta$ be the $n$--th root of unity $\,e^{2 \pi i/n}\,$ and 
$\,u',v'\,$ be the inverses of $\,u,v\,$ in $\,({\bf Z}/n{\bf Z})^*\,$. 
At each of their their fixed points, $g$ has the multiplier $\zeta^{u'}$ and 
$gh$ has the multiplier $\zeta^{v'}\,$. $\Box$
\end{lemma}

The parameters $u$ and $u'$ determine each other uniquely, 
and the same is true for $v$ and $v'\,$. Therefore, according to 
Theorem 1 the isomorphism class of $X(f;u,v,w)$ is uniquely determined 
by the pair of exponents $\,u',v' \bmod p^{e-f} \,$. Because under the  
action of all $\,\sigma \in {\rm Gal \,}\overline{{\bf Q}}/{\bf Q}(\eta)\,$ 
on the multipliers these residue classes remain unchanged, we have 
\[ X(f;u,v,w)^{\sigma} \; \cong \; X(f;u,v,w) \; . \]
By definition, the fixed field of all $\sigma$ satisfying this property 
is the {\em moduli field} $M$ of $X(f;u,v,w)\,$. In other words, we know 
now that $\,M \subseteq {\bf Q}(\eta) \,$. Since $X(f;u,v,w)$ is a quasiplatonic 
curve, we know by [W1, Remark 4] (for a more complete proof see [W2, Thm.~5]) 
that $X(f;u,v,w)$ can be defined over $M\,$. To finish the proof of 
Theorem 2, we may restrict 
all $\,\sigma \in {\rm Gal\,}\overline{{\bf Q}}/{\bf Q}\,$ to the cyclotomic field 
${\bf Q}(\eta)\,$, identify the Galois group $\,{\rm Gal \,}{\bf Q}(\eta) 
/{\bf Q}\,$ with $\,({\bf Z}/p^{e-f} {\bf Z})^*\,$ and recall that the 
action $\,\sigma(\eta) = \eta^{k'}\,$ extends to an action $\,\sigma( 
\zeta) = \zeta^{k'}\,$ if we replace $\,k'\in ({\bf Z}/p^{e-f} {\bf
  Z})^*\,$ with an integer representative, 
thus acting on the multipliers by 
\[ \zeta^{u'} \; \mapsto \; \zeta^{k'u'} \quad , \qquad \zeta^{v'} 
\; \mapsto \; \zeta^{k'v'} \;.\]
Taking a solution $\,k \in ({\bf Z}/p^{e-f}{\bf Z})^*\,$ of $\,k'k 
\equiv 1 \bmod p^{e-f}\,$, the Galois conjugation $\sigma$ acts on the 
triples $\,(u,v,w)\,$ (modulo $p^{e-f}$ and permutation) by 
\[ (u,v,w) \; \mapsto \; (ku,kv,kw) \; . \quad \Box \] 

To prove Theorem 3, we have to determine the moduli field $M$ more 
precisely, given a triple $\,(u,v,w) \bmod p^{e-f}\,$. In other words, 
we have to determine the subgroup of all $\,k \in ({\bf Z}/p^{e-f}{\bf Z})^*
\,$ with the property that modulo $p^{e-f}\,$, the triple $\,(ku,kv,kw)\,$ 
is just a permutation of $\,(u,v,w)\,$. This is possible only in the 
following situations. 
\begin{itemize}
\item $k=1\,$. If this trivial solution is the only one, we are in 
case 3 of Theorem 3.
\item $k$ of order $2\,$, i.e.~$\,k= -1\,$, hence in our normalisation 
for triples $\,(u,-u,0)\,$ with moduli field $\,M=R\,$ (case 1),
\item $k$ of order $3\,$, a solution of $\,1+k+k^2 \equiv 0 \bmod 
p^{e-f}\,$, and parameter triples $\,(u,ku,k^2u)\,$. Since we exclude 
the trivial situation $\,p^{e-f}=3 \,,\, k=1\,$, such curves 
exist if and only if $\,p \equiv 1 \bmod 3\,$ (case 2). 
$\Box$
\end{itemize}

{\em Examples.} For each $\,f<e\,$, there is one Galois orbit of curves 
$X(f;u,-u,0)$ defined over $R\,$, and for each $f$ there is one Galois orbit 
of curves $X(f;u,u,-2u)$ defined over the full cyclotomic field, 
treated in detail in [JStW] and already 
mentioned at the end of Section 2. The first example of case 2 occurs 
for $\,n=49 \,,\,p=7\,,\,f=1\,$, the Galois orbit consisting here of 
$\,X(1;1,2,4)\,,\,X(1;3,6,5)\,$, both being defined over $\,K={\bf Q}
(\sqrt{-7}) \subset {\bf Q}(e^{2\pi i/7})\,$ and hence conjugate under 
complex conjugation. 

{\em An alternative proof of the ``only if'' part of Theorem 1.}
Suppose that there is an isomorphism $\,F: X(f;u_1,v_1,w_1) \to
X(f;u_2,v_2,w_2)\,$. By permutation of the fixed points we may assume that
$\, u_1,v_1,u_2,v_2\,$ are coprime to $p\,$, and by composition with
automorphisms of $G_f$ we may assume that for both curves the group
action is normalised as in the beginning of this section,
i.e. that $g$ and $gh$ have fixed points of order $n\,$. Moreover, by
permutation of the fixed points and composition of $F$ with curve automorphisms
we may assume that $F$ maps fixed points $\,P_1,Q_1\,$ of $\,g, gh\,$ 
on the first curve to fixed points $\,P_2,Q_2\,$ of $g$ and $gh$ on the
second curve. Geometrically it is evident that for some integers $i$ and
$j$ 
\[ F \circ g \; = \; g^i \circ F \quad , \quad F \circ gh \; = \; (gh)^j
\circ F \;.\]
Therefore, $\,a \mapsto F \circ a \circ F^{-1}\,$ defines an
automorphism of $G_f$ sending $g$ and $gh$ to $g^i$ and
$(gh)^j\,$. According to [JN\v{S}, Prop.~16] this is possible only for $\,i,j
\equiv 1 \bmod p^{e-f}\,$, so we can assume even
\[ F \circ g \; = \; g \circ F \quad , \quad F \circ gh \; = \; (gh)
\circ F \;.\]
But then it is again geometrically evident that $g$ and $gh$ have the same 
multipliers at $\,P_1,Q_1\,$ as at $\,P_2,Q_2\,$, so $\,(u_1,v_1,w_1) =
(u_2,v_2,w_2)\,$ follows from Lemma 3. $\Box$


\section{Equations}

\begin{thm}
Let $\,n=p^e\,$ be an odd prime power. Suppose that $\,2f\ge e\,$ 
and (without loss of generality) that $u,v$ are not divisible by $p\,$. Then 
an affine model of $\,X(f;u,v,w)\,$ in $ {\bf C}^4$ is given by 
the equations 
\begin{eqnarray}
y^n & = & \beta^{u}(1-\beta)^{v} \\
x^{p^{e-f}} & = & 1-\beta \\
z^{p^f} & = & x^{-r} \prod_{m=0}^{p^{e-f}-1} (x-\eta^m)^{am} \;,
\end{eqnarray}
where $\,a:= p^{2f-e} \,$ and the exponent $\,r:=(q^{p^{e-f}}-1)/p^e\,$ 
is an integer coprime to $p\,$.
\end{thm} 

{\em Proof.} $\,H := \langle h \rangle\,$ is the normal subgroup of  
$G_f$ whose quotient group $\,C_n \cong G_f/H\,$ is generated by $gH\,$. 
Its preimage under the epimorphism $\,\Delta \to G_f\,$ is a Fuchsian 
normal subgroup $\Gamma_H$ of $\Delta$ with quotient curve $X_H\,$. 
This curve is a cyclic cover of the projective line 
$\,{\bf P}^1({\bf C}) \cong \Delta \backslash {\bf H}\,$ of degree $n$ 
whose function field is therefore a cyclic extension of a rational 
function field ${\bf C}(\beta)$ determined by an equation 
$\,y^n = f(\beta) \in {\bf C}(\beta)\,$. Since the covering is ramified
only 
over $\,\beta = 0,1,\infty \,$, we may suppose that $X_H$ is a Weil
curve 
\[ y^n = \beta^c (1-\beta)^d \]
and that the covering map $\,X_H \to {\bf P}^1({\bf C}) \,: 
\,(y,\beta) \mapsto \beta\,$ is a Bely\u{\i} function on
$X_H\,$ whose branches are in bijective correspondence with the cosets 
$\,Hg, Hg^2,\ldots, Hg^{n}=H\,$ of $H\,$. Without loss of generality, 
we may assume that $\,\beta = 0\,$ at the image points of the $\Delta$--orbit 
$\Delta z_0$ of the fixed point of $\gamma_0$  under the quotient 
map $\,{\bf H} \to X_H\,$, 
and $\,\beta = 1\,$ at those coming from the 
$\Delta$--orbit $\Delta z_1$ of the fixed point of $\gamma_1\,$. In fact, 
$\beta$ is induced by a $\Delta$--automorphic function on ${\bf H}$ 
mapping 
the two (open) triangles forming the fundamental domain for $\Delta$ 
conformally onto the upper and the lower half plane. 

Let $\beta^{\frac{1}{n}}$ be a branch of an $n$--th root of $\beta$ 
multiplied by the factor $\zeta =  e^{2 \pi i/n}$  by counterclockwise
continuation around $0\,$, and similarly $(1-\beta)^{\frac{1}{n}}$ 
by continuation around $1\,$. The exponents $c$ and $d$ have to be
chosen so that this description of the covering surface corresponds to the
labelling by the cosets chosen above. Since the canonical epimorphism $\,\Delta
\to \Delta/\Gamma_H \cong G_f/H \cong C_n\,$ is determined by 
\[ \gamma_0 \mapsto (Hg)^u = Hg^u\quad , \qquad \gamma_1 \mapsto (Hg)^v =
Hg^v \;, \] 
counterclockwise continuation along paths around $0$ and $1$ has to give
cycles of branches
$\,(H,Hg^u,Hg^{2u}, \ldots)\,$ and $\,(H, Hg^v, Hg^{2v},
\ldots)\,$, respectively. Therefore the passage from the branch $H$ to
the branch $Hg^{ju}$ corresponds to the factor $\zeta^j$ for
$\beta^{\frac{1}{n}}\,$, and  the passage from the branch $H$ to
the branch $Hg^{jv}$ corresponds to the factor $\zeta^j$ for
$(1-\beta)^{\frac{1}{n}}\,$.  
We are free to choose the exponent $c$ to be any number coprime to
$n\,$, so we can take $\,c = u\,$. Then the passage from  the branch $H$ to
the branch $Hg^{ju}$ corresponds to the factor $\zeta^{ju}$ for
$\beta^{\frac{u}{n}}$ and similarly the passage from the branch $H$ to
the branch $Hg^{jv}$ corresponds to the factor $\zeta^{jv}$ for
$(1-\beta)^{v\over n}\,$. Hence $\,d=v\,$ is consistent with the
choice $\,c=u\,$, leading to the same factor for the continuation of $\,y =
\sqrt[n]{\beta^u(1-\beta)^v} \,$ on any path in the $\beta$--plane
avoiding $0$ and $1\,$.

Up to this point, we have not used the hypothesis that $\,2f \ge e\,$: the Weil
curve (1) is in fact always a quotient of $X(f;u,v,w)\,$. The rest of
the proof does not differ from [JStW, Lemmata 9 to 11], so it may be
sufficient to sketch the remaining part. Equation (2)
describes a cyclic degree $p^{e-f}$ extension of the base function field
${\bf C}(\beta)\,$. 
The extended field is the fixed field of the subgroup of $G_f$ generated by
$g$ and $h^{p^{e-f}}$ and is of genus $0\,$, and equation (3) describes in
turn a cyclic extension of this field, i.e.~the fixed field of the subgroup
$\langle g \rangle\,$. Here, the hypothesis $\,2f \ge e\,$ plays an
essential role. Since $H$ and $\langle g \rangle$ generate $G_f$ and
have trivial intersection, all three equations together describe the
curve. Clearly, $\beta$ can be eliminated by (2). $\Box$

Equations (1) and (2) are defined over the rationals and hence describe 
quotient curves for all curves of one Galois orbit. In the cases 1 and 2
of Theorem 3, equation (3) is not defined over the minimal possible 
field of definition since --- according to Theorem 3 --- the curves can
then be defined 
over proper subfields of ${\bf Q}(\eta)\,$. In principle, such 
a model defined over the field of moduli can be found by 
a different choice of the coordinates better reflecting the symmetry 
between the critical values of the Bely\u{\i} function, see [StW1, Remark 1] 
or the proof of [W2, Thm.~5]. 

\vspace{1cm}


{\bf References}

[BCIR] M.~Bauer, A.~Coste, Cl.~Itzykson, P.~Ruelle, Comments on the
links between ${\rm su}(3)$ modular invariants, simple factors in the
Jacobians of Fermat curves, and rational triangular billards,
J. Geom. Phys. 22 (2) (1997), 134--189.

[Bel] G.V.~Bely\u{\i}, On Galois extensions of a maximal cyclotomic
field, Math.~USSR Izvestija 14 (1980) 247--256.

[BCC] E.~Bujalance, F.J.~Cirre, M.~Conder, On Extendability of Group
Actions on Compact Riemann Surfaces, Transactions of the AMS 355 No.4
(2002) 1537--1557.

[CIW] P.~Cohen, Cl.~Itzykson, J.~Wolfart, Fuchsian triangle groups and
Grothendieck dessins: variations on a theme of Belyi,
Commun.~Math.~Phys. 163 (1994) 605--627.

[G] A.~Grothendieck, Esquisse d'un Programme, pp.~5--84 in {\em Geometric
Galois Actions 1. Around Grothendieck's Esquisse d'un Programme},
ed. P.~Lochak, L.~Schneps, London Math.~Soc.~Lecture Note Ser. 242, Cambridge
University Press, Cambridge, 1997.

[GW] E.~Girondo, J.~Wolfart, Conjugators of Fuchsian groups and
quasiplatonic surfaces, Quarterly J.~Math. 56 (2005) 525--540 

[JN\v{S}] G.A.~Jones, R.~Nedela, M.~\v{S}koviera, Regular embeddings of 
$K_{n,n}$ where $n$ is an odd prime power, European J.~Combinatorics, to
appear.

[JS] G.A.~Jones, D.~Singerman, Belyi functions, hypermaps and Galois 
groups, Bull.~London Math.~Soc. 28 (1996) 561--590. 

[JStW] G.A.~Jones, M.~Streit, J.~Wolfart, Galois action on families of 
generalised Fermat curves, J. of Algebra 307 (2007) 829--840. 

[Si] D.~Singerman, Finitely maximal Fuchsian groups, 
J.~London Math.~Soc.(2) 6 (1978) 29--38.

[St] M.~Streit, Field of definition and Galois orbits for the 
Macbeath--Hurwitz curves, Arch.~Math. 74 (2000) 342--349. 

[StW1] M.~Streit, J.~Wolfart, Characters and Galois invariants of regular 
dessins, Revista Mat.~Complutense 13 (2000) 49--81. 

[StW2] M.~Streit, J.~Wolfart, Cyclic Projective Planes and Wada Dessins,
Documenta Mathematica 6 (2001) 39--68.

[VS] V.A.~Voevodsky, G.~Shabat, Equilateral triangulations of Riemann
surfaces and curves over algebraic number fields, Soviet Math.~Dokl. 39
(1989) 38--41.

[W1] J.~Wolfart, The `Obvious' part of Belyi's Theorem and Riemann 
surfaces with many automorphisms, pp.~97--112 in Geometric Galois 
Actions 1, ed.~L.~Schneps and P.~Lochak, London Math.~Soc.~Lecture Note Ser. 242, 
Cambridge University Press, Cambridge, 1997.

[W2] J.~Wolfart, ABC for polynomials, dessins d'enfants, and
uniformization --- a survey, pp. 313--345 in {\em Elementare und
  Analytische Zahlentheorie (Tagungsband),} Proceedings ELAZ--Conference
  May 24--28, 2004
(ed. W.~Schwarz, J.~Steuding), Steiner Verlag Stuttgart 2006
({\tt http://www.math.uni-frankfurt.de/$\sim$wolfart/}). 

\end{document}